\documentclass[10pt]{amsart}
\usepackage{a4}
\usepackage{amssymb}
\newtheorem{thm}{Theorem}

\newtheorem{prop}[thm]{Proposition}
\newtheorem{defn}[thm]{Definition}

\newtheorem{rem}[thm]{Remark}
\newtheorem*{tthm}{Theorem}

\numberwithin{thm}{section}
\numberwithin{equation}{section}

\newcommand{\norm}[1]{\left\Vert#1\right\Vert}
\newcommand{\abs}[1]{\left\vert#1\right\vert}

\newcommand{\F}{\mathcal{F}}
\newcommand{\G}{\mathcal{G}}
\newcommand{\Hi}{\mathcal{H}}

\newcommand{\M}{\mathcal{M}}
\newcommand{\N}{\mathcal{N}}

\newcommand{\U}{\mathcal{U}}

\newcommand{\W}{\mathcal{W}}
\newcommand{\WT}{\widetilde{\mathcal{W}}}

\newcommand{\E}{\mathbb{E}}
\newcommand{\n}{\mathbb{N}}

\newcommand{\Om}{\Omega}
\newcommand{\om}{\omega}

\begin{document}
\title[Unconditionality]
{Unconditionality with respect to orthonormal systems in noncommutative $L_2$ spaces}
\author{Hun Hee Lee}

\address{Department of Mathematical Sciences, Seoul National University
          San56-1 Shinrim-dong Kwanak-gu Seoul 151-747, Republic of Korea}

\begin{abstract}
Orthonormal systems in commutative $L_2$ spaces can be used to classify Banach spaces. 
When the system is complete and satisfies certain norm condition the unconditionality with respect to the system characterizes Hilbert spaces. 
As a noncommutative analogue we introduce the notion of unconditionality of operator spaces with respect to orthonormal systems 
in noncommutative $L_2$ spaces and show that the unconditionality characterizes operator Hilbert spaces 
when the system is complete and satisfy certain norm condition. The proof of the main result heavily depends on free probabilistic tools such as 
contraction principle for $*$-free Haar unitaries, comparision of averages with respect to $*$-free Haar unitaries and $*$-free circular elements 
and $K$-covexity, type 2 and cotype 2 with respect to $*$-free circular elements.

\end{abstract}

\maketitle

\section{Introduction}
In Banach space theory orthonormal systems in classical $L_2$ spaces have been frequently used for many purposes. 
For an orthonormal system $\Phi = (\phi_i)_{i\geq 1}$ in $L_2(M, \mu)$, where $(M, \mu)$ is a measure space, we consider
so called ``$\Phi$-average" of a sequence $(x_i)^n_{i=1}$ in a Banach space $X$ as follows. 
$$\norm{\sum^n_{i=1}\phi_i \otimes x_i}_{L_2(\mu;X)} = 
\Big[\int_M \norm{\sum^n_{i=1}\phi_i(t)x_i}_X^2d\mu(t)\Big]^{\frac{1}{2}}.$$
Many properties of Banach space are defined using this ``averages". 
For example, type 2 and cotype 2 conditions are defined by comparing averages with respect to 
rademacher systems $(\varepsilon_i)_{i\geq 1}$ in $L_2[0,1]$ and unit vector systems $(e_i)_{i\geq 1}$ in $\ell_2$, 
where $\varepsilon_i(t) = \text{sign}(\sin(2^i\pi t))$, $t\in [0,1]$ and $i=1,2,\cdots$.

Thus, it is natural to expect orthonormal systems in noncommutative $L_2$ spaces play 
a similar role in operator space theory. In this paper we focus on an operator space version of the result in \cite{DJ92}.
In \cite{DJ92} another example of using averages, namely, the notion of $\Phi$-unconditionality was considered.
We say that a Banach space $X$ is $\Phi$-unconditional if $$\norm{\sum^n_{i=1}\epsilon_i\phi_i \otimes x_i}_{L_2(\mu;X)} 
\sim \norm{\sum^n_{i=1}\phi_i \otimes x_i}_{L_2(\mu;X)}$$ for any
$n\in \n$, $(x_i)^n_{i=1}\subseteq X$ and any $\epsilon_i \in \{\pm 1\}$, where `$\sim$' implies equivalent allowing constants.
The authors showed that if $(\phi_i)_i$ is complete and $$\norm{ \sup_i\abs{\phi_i}}_{L_2(\mu)} = \Big[\int_M \norm{\sup_i\abs{\phi}_i(t)}^2d\mu(t)\Big]^{\frac{1}{2}} < \infty,$$
then a Banach space is $\Phi$-unconditional if and only if it is isomorphic to a Hilbert space.
 
In order to consider an operator space analogue of the above result we first need to define 
the right notion of unconditionality. Let $\M$ be a von Neumann algebra and $\Phi=(\phi_i)_{i\geq 1}$ be an orthonormal system in $L_2(\M)$. 
We will define that an operator space $E$ is $\Phi$-unconditional if 
$$\norm{\sum^n_{i=1} \phi_i \otimes u_i \otimes x_i}_{L_2(\M \overline{\otimes}M_m ;E)} 
\sim \norm{\sum^n_{i=1} \phi_i \otimes x_i}_{L_2(\M;E)}$$ for any $n,m \in \n$, $(x_i)^n_{i=1}\subseteq E$ and any unitaries $u_i \in M_m$. 
Actually, we need some restrictions on $\M$ and $E$ to make $L_2(\M; E)$ sense.
Recall that a $C^*$-algebra $A$ has $WEP$ (Lance's weak expectation property) if the inclusion map $i_A : A \hookrightarrow A^{**}$ 
factors completely positively and completely contractively through $B(H)$ for some Hilbert space $H$.
A $C^*$-algebra $B$ has $QWEP$ if there is a $WEP$ $C^*$-algebra $A$ and two sided ideal $I\subseteq A$ such that $B \cong A/I$.
Let $C^*(\mathbb{F}_{\infty})$ be the full group $C^*$-algebra of the free group with infinite number of generators.
Then we say that an operator space $E$ is locally-$C^*(\mathbb{F}_{\infty})$ (chapter 21. of \cite{P03} and \cite{Har98}) if 
$$d_f(E) := \sup_{F\subseteq E, \text{finite dimensional}} \inf\{ d_{cb}(F,G) 
: G \subseteq C^*(\mathbb{F}_{\infty}) \} < \infty.$$
If $\M$ has $QWEP$ and $E$ is locally-$C^*(\mathbb{F}_{\infty})$, then the above $L_2(\M; E)$ is well-defined. (See section \ref{sec-prelim} for the details.)

In order to impose a norm condition on $\Phi$ we need another concept. Recall that (see \cite{J02} for the details) 
for a sequence $(a_i)_{i\geq 1}$ in $L_p(\M)$ $(1 \leq p < \infty)$ we define 
$$\norm{\sup_i \abs{a_i}}_p := \inf_{a_i = ay_i b} \norm{a}_{L_{2p}(\M)}\norm{b}_{L_{2p}(\M)}\sup_i \norm{y_i}_{\M}.$$
Note that the notation $\sup_i a_i$ was used in \cite{J02} instead of $\sup_i \abs{a_i}$ but we will use the latter for the compatibility with commutative cases. 

Then, the main result is as follows. 
\begin{tthm}
Let $\M$ be a von Neumann algebra with $QWEP$ which is not subhomogeneous, and $\Phi=(\phi_i)_{i\geq 1}$ be
a complete orthonormal system in $L_2(\M)$ with $$\norm{\sup_i \abs{\phi_i}}_2 < \infty.$$
Then, a locally-$C^*(\mathbb{F}_{\infty})$ operator space $E$ is $\Phi$-unconditional 
if and only if $E$ is completely isomorphic to an operator Hilbert space.
\end{tthm}
We exclude the case $\M$ is subhomogeneous since the unconditionality in this case does not use the whole operator space structure of $E$.  

The essential tools for the proof of Banach space case is concerned with averages 
with respect to the Rademacher system and the standard gaussian variables. 
Thus, it is also natural to expect to their noncommutative analogues, 
$*$-free Haar unitaries and $*$-free circular elements would paly a similar role in oeprator space case.

This paper is organized as follows. In section \ref{sec-prelim} we collect some basic materials
about vector valued noncommuatative $L_p$ spaces with respect to $QWEP$ von Neumann algebras,
$*$-free Haar unitaries and $*$-free circular elements. In section \ref{sec-tools} we develop
the following essential tools for the proof of the main results; contraction principle for $*$-free Haar unitaries,
comparision of averages with respect to $*$-free Haar unitaries and $*$-free circular elements and
$K$-covexity, type 2 and cotype 2 with respect to $*$-free circular elements. 
In the last section we define the unconditionality of operator spaces with respect to 
orthonormal systems in noncommutative $L_2$ spaces and prove the main result.

Throughout this paper, we assume that the reader is familiar with the general results of 
operator spaces (\cite{ER00,P03}), operator algebras (\cite{Ta79}), free probability (\cite{HiPe00, VDN92}) and completely $p$-summing maps (\cite{P98}).
For a index set $I$ we denote the operator Hilbert space on $\ell_2(I)$ by $OH(I)$. When $I = \{1, \cdots, n\}$, $n\in \n$ we simply write $OH_n$. 

\section{Preliminaries}\label{sec-prelim}

\subsection{Vector valued noncommuatative $L_p$ spaces 
with respect to $QWEP$ von Neumann algebra}\label{subsec-LpQWEP}

The theory of vector valued noncommuatative $L_p$ spaces was initiated by G. Pisier (\cite{P98})
for the case that the underlying von Neumann algebra is injective and semifinite. For an operator space $E$ and an injective and semifinite 
von Neumann algebra $M$ we define $$L_p(M; E):= [M\otimes_{\min} E, L_1(M)\widehat{\otimes} E]_{\frac{1}{p}},$$ where $\otimes_{\min}$ and $\widehat{\otimes}$ imply injective and projective tensor product of operator spaces and $[\cdot, \cdot]_{\theta}$ 
implies complex interpolation of operator spaces. 
Recently, M. Junge (\cite{J-QWEP}) extended this theory to the case that the underlying von Neumann algebra satisfies $QWEP$ using the following characterization of $QWEP$ von Neumann algebras.

\begin{prop}\label{prop-QWEP-Equi}
A von Neumann algebra $\M$ is $QWEP$ if and only if there are a normal $*$-isomorphism 
$$\pi : \M \rightarrow \M_{\U} := \Big(\prod_{i, \U}L_1(M_i)\Big)^*$$ for some injective and semifinite von neumann algebras $M_i$ 
and some free ultrafilter $\U$ on an index set $I$ and a normal conditional expectation
$$\mathcal{E} : \M_{\U} \rightarrow \pi(\M).$$
\end{prop}

Let $\M$ be a von Neumann algebra with $QWEP$ and $\alpha = (\pi, \mathcal{E})$ as above.
Then for $1\leq p \leq \infty$ there are complete isometries 
$$\pi_p : L_p(\M)\rightarrow L_p(\M_\U) = \prod_{\U}L_p(M_i)$$
and complete contractions $$\mathcal{E}_p : L_p(\M_\U) \rightarrow L_p(\pi(\M))$$ 
induced from $\pi$ and $\mathcal{E}$, respectively.
The followings are basic properties of $L_p(\M, \alpha ;E)$ from \cite{J-QWEP} which we will need in the sequel.

\begin{prop}\label{prop-LpQWEP-Basic}
Let $\M$ be a von Neumann algebras with $QWEP$, $\alpha = (\pi, \mathcal{E})$ as above, 
$E_1$ and $E_2$ be operator spaces, $1\leq p < \infty$ and $\frac{1}{p} + \frac{1}{p'} =1$.
\begin{itemize}
\item[(1)]
If $E_1 \subseteq E_2$ is a completely isometric embedding, then we have
$$L_p(\M, \alpha ;E_1)\subseteq L_p(\M, \alpha ;E_2)$$ is completely isometric.

\item[(2)]
If $T : E_1 \rightarrow E_2$ is a completely bounded map, then we have
$$\norm{I_{L_p(\M)}\otimes T : L_p(\M, \alpha ;E_1)\rightarrow L_p(\M, \alpha ;E_2)} \leq \norm{T}_{cb}.$$

\end{itemize}

\end{prop}

The definition of $L_p(\M, \alpha ;E)$ depends on the choice of $\alpha$. 
However, it can be simplified if we restrict our attention to the case $E \subseteq C^*(\mathbb{F}_{\infty})$.
Indeed, we have $$\norm{x}_{L_p(\M,\alpha;E)} = \inf_{x = a\cdot y\cdot b} 
\norm{a}_{L_{2p}(\M)} \norm{y}_{\M \otimes_{\max}C^*(\mathbb{F}_{\infty})}\norm{b}_{L_{2p}(\M)}$$
for $x \in L_p(\M)\otimes E$ and $1\leq p < \infty$, where $\otimes_{\max}$ implies maximal $C^*$-tensor product. 
In particular, $L_p(\M,\alpha;E)$ does not depend on the choice of $\alpha$. In this case, we use simpler notation $L_p(\M;E)$.

When $E$ is locally-$C^*(\mathbb{F}_{\infty})$ by Proposition \ref{prop-LpQWEP-Basic} 
$L_p(\M, \alpha ;E)$'s are all equivalent allowing the constant $d_f(E)^2$. 
Thus, we can still say that $L_p(\M,\alpha;E)$ does not depend on the choice of $\alpha$, and thus we also use the notation $L_p(\M;E)$. 
Note that exactness implies ``locally-$C^*(\mathbb{F}_{\infty})$", and ``locally-$C^*(\mathbb{F}_{\infty})$" is stable under duality 
and complex interpolation (\cite{Har98}), so that $S_p$ and $L_p(\mu)$ ($1\leq p \leq \infty$) for some measure $\mu$ are all locally-$C^*(\mathbb{F}_{\infty})$.

Furthermore, the condition $E \subseteq C^*(\mathbb{F}_{\infty})$ guarantees the $E$-valued extension property 
of completely positive maps between noncommutative $L_p$-spaces. (\cite{J-QWEP}) 
Recall that, for von Neumann algebras $\M_1$ and $\M_2$, $S : L_p(\M_1) \rightarrow L_p(\M_2)$ is called completely positive
if for each $n \in \n$, $I_{M_n}\otimes S$ maps the positive cone $L_p(M_n \overline{\otimes} \M_1)^+$ 
into the positive cone $L_p(M_n \overline{\otimes} \M_2)^+$.

\begin{prop}\label{prop-CPmap}
Let $\M_1$ and $\M_2$ be von Neumann algebras with $QWEP$, 
$E \subseteq C^*(\mathbb{F}_{\infty})$ and $1\leq p < \infty$.
If $S : L_p(\M_1) \rightarrow L_p(\M_2)$ is a completely positive map, then we have
$$\norm{S \otimes I_E : L_p(\M_1 ;E)\rightarrow L_p(\M_2 ;E)} \leq \norm{T}.$$

\end{prop}

We end this subsection with the following duality result.

\begin{prop}\label{prop-duality}
Let $\M$ be a von Neumann algebra with $QWEP$, $\alpha = (\pi, \mathcal{E})$ as above, 
$E$ be a locally-$C^*(\mathbb{F}_{\infty})$ operator space, $1< p < \infty$ and $\frac{1}{p} + \frac{1}{p'} = 1$.
Then, we have $$L_{p'}(\M,\alpha ; E^*) \hookrightarrow L_p(\M, \alpha ; E)^*$$ completely isomorphically. 
More precisely, we have $$\norm{\xi}_{L_{p}(\M, \alpha; E)^*} 
\leq \norm{\xi}_{L_{p'}(\M, \alpha; E^*)} \leq d_f(E)\norm{\xi}_{L_{p}(\M, \alpha; E)^*}$$
for any $\xi \in L_{p'}(\M, \alpha; E^*)$.
\end{prop}
\begin{proof}
We consider $\xi \in L_{p'}(\M)\otimes E^*$ with $\pi_{p'}\otimes I_E(\xi) = (\xi_i)_\U \in \prod_\U L_{p'}(M_i; E^*)$, 
then for $x \in L_p(\M)\otimes E$ with $\pi_{p}\otimes I_E(x) = (x_i)_\U \in \prod_\U L_p(M_i; E)$ we have
\begin{align*}
\begin{split}
\abs{\left\langle \xi, x \right\rangle} & = \abs{\left\langle \pi_{p'}(\xi), \pi_{p}(x) \right\rangle}
= \abs{\lim_\U \left\langle \xi_i, x_i \right\rangle}\\
& \leq \lim_\U \norm{\xi_i}_{L_{p'}(M_i; E)} \lim_\U \norm{x_i}_{L_p(M_i; E)}\\
& = \norm{\xi}_{L_{p'}(\M, \alpha; E^*)} \norm{x}_{L_{p}(\M, \alpha; E)},
\end{split}
\end{align*}
which implies $\norm{\xi}_{L_{p}(\M, \alpha; E)^*} \leq \norm{\xi}_{L_{p'}(\M, \alpha; E^*)}$.

For the converse inequality we choose $y_i \in L_p(M_i; E)$ with 
$$\left\langle \xi_i, y_i \right\rangle = \norm{\xi_i}_{L_{p'}(M_i; E^*)}\,\, \text{and}\,\,
\norm{y_i}_{L_p(N_i)} \leq 1 + \frac{1}{i}.$$
Let $\pi_{p'} : L_{p'}(\M)\rightarrow \prod_{\U}L_{p'}(M_i)$ be the complete isometry deriveded from $\pi$, 
which is clearly completely positive. Then, its adjoint 
$\pi^*_{p'} : \prod_{\U}L_{p}(M_i)\rightarrow L_{p}(\M)$ is also completely positive and
\begin{align*}
\begin{split}
\abs{\left\langle \xi, \pi^*_{p'}\otimes I_E [(y_i)_\U] \right\rangle}
& = \abs{\left\langle \pi_{p'}\otimes I_E (\xi), (y_i)_\U \right\rangle}\\
& = \abs{\lim_\U \left\langle \xi_i, y_i \right\rangle} = \lim_\U \norm{\xi_i}_{L_{p'}(M_i; E^*)}\\
& = \norm{\xi}_{L_{p'}(M_i ;E^*)}.
\end{split}
\end{align*}
On the other hand we have by Proposition \ref{prop-CPmap} 
$$\norm{\pi^*_{p'}\otimes I_E[(y_m)_\U]}_{L_{p}(\M; E)}
\leq d_f(E)\norm{\pi^*_{p'}}\lim_\U\norm{y_m}_{L_p(M_i; E)} \leq d_f(E),$$
which implies $\norm{\xi}_{L_{p'}(N;E^*)}\leq d_f(E)\norm{\xi}_{L_{p}(N;E)^*}$.

By repeating the same argument for $S^n_p(E)$, $n \geq 1$ instead of $E$ we get the desired complete isometry. 

\end{proof}

\subsection{$*$-free Haar unitaries and $*$-free circular elements}\label{subsec-RandomMatrix}

In this subsection we consider specific $\alpha = (\pi, \mathcal{E})$'s for free group von Nemann algebra
obtained by using random matrix models for $*$-free Haar unitaries and $*$-free circular elements.

Let $\mathbb{F}_\infty$ be the free group with generators $(g_i)_{i\geq 1}$ and
$\lambda(g_i)$ be the left translation by $g_i$ in $\ell_2(\mathbb{F}_{\infty})$ for $i\geq 1$.

Let $\Hi$ be a Hilbert space with Hilbert space basis $(e_i)_{i\geq 1} \cup (f_i)_{i\geq 1} \subseteq \Hi$. 
Now we consider the full Fock space $$\F(\Hi) := \mathbb{C}\Omega \oplus_{n \geq 1}\Hi^{\otimes n},$$
and let $\WT_i := \ell(e_i) + \ell^*(f_i)$, where $\ell(f) \in B(\F(\Hi))$, $f\in \Hi$ is the left creation operator
defined by $$\ell(f)(\Omega) := f$$ and $$\ell(f) (f_1\otimes \cdots \otimes f_n) := f\otimes f_1\otimes \cdots \otimes f_n$$ 
for $n\geq 1$ and $f_1, \cdots, f_n \in \Hi$, and $\ell^*(f)\in B(\F(\Hi))$ is the adjoint of $\ell(f)$.

$(\lambda(g_i))_{i\geq 1}$ and $(\WT_i)_{i\geq 1}$ are typical examples of 
$*$-free Haar unitaries and $*$-free circular elements, respectively, 
and they have the random matrix models as follows.
Let $(\Om, P)$ and $(\Om', P')$ be probability spaces, $m \in \n$ and 
$\U(m)$ be the compact group of $m \times m$ unitary matrices with the normalized Haar measure $\gamma_n$. 
Now we consider the standard unitary random matrix $$U^m : (\Om, P) \rightarrow \U(m)$$ 
with distribution $\gamma_n$ and the standard gaussian random matrix 
$$G^m : (\Om', P') \rightarrow M_m,\,\, \text{where}\,\, G^m = \Big(\frac{1}{\sqrt{m}}g_{ij}\Big)^m_{i,j=1}$$ 
and $g_{ij}$'s are i.i.d. complex valued standard gaussian random variables.
 
Then $U^m$ and $G^m$ are noncommutative random variables in $$(L_{\infty}(\Om ;M_m), \tau_m)\,\, 
\text{and}\,\, (\cap_{1\leq p < \infty}L_p(\Om' ;M_m), \tau'_m),$$ respectively, where 
$$\tau_m(x) = \int_\Om \frac{1}{m}\text{tr}(x(\om))dP(\om)\,\, \text{and}\,\, 
\tau'_m(y) = \int_{\Om'} \frac{1}{m}\text{tr}(y(\om'))dP'(\om')$$ for $x : \Om \rightarrow M_m$ and $y : \Om' \rightarrow M_m$,
and it is well known (by Voiculescu) that $(U^m_i)_{i\geq 1}$ (resp. $(G^m_i)_{i\geq 1}$), independent copies of $U^m$ (resp. $G^m$), 
converges in distribution to $(\lambda(g_i))_{i\geq 1}$ (resp. $(\WT_i)_{i\geq 1}$) as $m$ goes to infinity.

Now we consider a free ultrafilter $\U$ in $\n$ and denote $\N_m := L_{\infty}(\Om ;M_m)$. 
Since $\Big(\prod_{\U}L_1(\N_m)\Big)^*$ coincide 
with the ultraproduct of $(\N_m)_{m\geq 1}$ in the sense of finite von Neumann algebras 
the above convergence implies that $$((U^m_i)_{\U})_{i\geq 1}\subseteq \Big(\prod_{\U}L_1(\N_m)\Big)^*$$
has the same $*$-distribution with $(\lambda(g_i))_{i\geq 1}$. Thus, we get a normal $*$-isomorphism
$$\pi_U : VN(\mathbb{F}_\infty) \rightarrow \Big(\prod_{\U}L_1(\N_m)\Big)^*$$ extending the natural map 
$$P(\lambda(g_1), \lambda(g_2), \cdots) \rightarrow P((U^m_1)_{\U}, (U^m_2)_{\U}, \cdots)$$ 
for any noncommutative polynomial $P$. (See Lemma 1 of \cite{DH03} for example) 
Moreover, since $\pi_U(VN(\mathbb{F}_\infty))\subseteq \Big(\prod_{\U}L_1(\N_m)\Big)^*$ are both finite von neumann algebras
with respect to the same trace we have the natural normal conditional expectation 
$\mathcal{E}_U :\Big(\prod_{\U}L_1(\N_m)\Big)^* \rightarrow \pi_U(VN(\mathbb{F}_\infty))$,
which is the adjoint of the inclusion $\pi_U(VN(\mathbb{F}_\infty))\hookrightarrow \Big(\prod_{\U}L_1(\N_m)\Big)^*$.

For the gaussian case we need to truncate since gaussian variables are not bounded.
Let $$\widetilde{G}^m : (\Om', P') \rightarrow M_m,\,\, \text{where}\,\, 
\widetilde{G}^m = \Big(\frac{1}{\sqrt{m}}1_{\sum^m_{i,j=1}\abs{g_{ij}}^2 \leq m^{2m}}\cdot g_{ij}\Big)^m_{i,j=1}.$$
Then for $N_m := L_{\infty}(\Om' ;M_m)$ we have $\widetilde{G}^m \in (N_m, \tau'_m)$, and 
$(\widetilde{G}^m_i)_{i\geq 1}$, independent copies of $\widetilde{G}^m$, 
has the same asymptotic $*$-distribution as $(G^m_i)_{i\geq 1}$.
Indeed, for fixed $k\in \n$ and $$\widehat{G}^m = G^m - \widetilde{G}^m = 
\Big(\frac{1}{\sqrt{m}} \widehat{g}_{ij}\Big)^m_{i,j=1},$$ we have 
$$\tau'_m((\widehat{G}^m)^k) = \frac{1}{m^{\frac{k}{2}+1}}\sum_{1\leq i_1,\cdots, i_k \leq m} 
\E(\widehat{g}_{i_1 i_2}\widehat{g}_{i_2 i_3}\cdots \widehat{g}_{i_m i_1}),$$ 
where $\E$ implies the expectation with respect to $(\Om', P')$. 
Since the cube $$[-m^{m-\frac{1}{2}}, m^{m-\frac{1}{2}}]^m \subseteq \mathbb{R}^m$$ is contained in 
the ball centered at $0$ with radius $m^m$ we have
\begin{align*}
\begin{split}
\abs{\E(\widehat{g}_{i_1 i_2}\widehat{g}_{i_2 i_3}\cdots \widehat{g}_{i_m i_1})}
& \leq \E(1_{\abs{g_{ij}} > m^{m-\frac{1}{2}}, \,\, 1\leq i,j \leq m}\cdot \abs{g_{i_1 i_2}g_{i_2 i_3}\cdots g_{i_m i_1}})\\
& \leq \E(1_{\abs{g_{i_1 i_2}} > m^{m-\frac{1}{2}}}\cdot \abs{g_{i_1 i_2}}
\cdots 1_{\abs{g_{i_m i_1}} > m^{m-\frac{1}{2}}}\cdot \abs{g_{i_m i_1}})\\
& \leq \Big[\E(1_{\abs{g_{i_1 i_2}} > m^{m-\frac{1}{2}}}\cdot \abs{g_{i_1 i_2}}^k)\Big]^{\frac{1}{k}}\cdots\\
& \,\,\,\,\,\,\,\,\Big[\E(1_{\abs{g_{i_m i_1}} > m^{m-\frac{1}{2}}}\cdot \abs{g_{i_m i_1}}^k)\Big]^{\frac{1}{k}}\\
& = \E(1_{\abs{g_{11}} > m^{m-\frac{1}{2}}}\cdot \abs{g_{11}}^k)
\leq C e^{-m^{m-\frac{1}{2}}}
\end{split}
\end{align*}
for some constant $C>0$ independent of $m$. Thus, we have 
$$\abs{\tau'_m((\widehat{G}^m)^k)} \leq C\frac{m^{m - \frac{k}{2} - 1}}{e^{m^{m-\frac{1}{2}}}} \rightarrow 0\,\,
\text{as}\,\, m \rightarrow \infty,$$ which implies $\widehat{G}^m$ converges to $0$ in $*$-distribution, 
and consequently $\widetilde{G}^m$ converges to $G^m$ in $*$-distribution.
On the other hand since $\widetilde{G}^m$ is bi-unitarily invariant 
(i.e. $V_1\widetilde{G}^m V_2$ ahs the same $*$-distribution as $\widetilde{G}^m$ 
for every $m \times m$ unitary matrices $V_1$ and $V_2$)  $(\widetilde{G}^m_i)_{i\geq 1}$ is asymptotically $*$-free 
(Theorem 4.3.11 of \cite{HiPe00}), so that we get the desired asymptotic $*$-distribution of $(\widetilde{G}^m_i)_{i\geq 1}$.

Then we get a normal $*$-isomorphism
$$\pi_G : \{\WT_i : i\geq 1\}'' \rightarrow \Big(\prod_{\U}L_1(N_m)\Big)^*$$ extending the natural map 
$$P(\WT_1, \WT2, \cdots) \rightarrow P((\widetilde{G}^m_1)_{\U}, (\widetilde{G}^m_2)_{\U}, \cdots)$$ 
for any noncommutative polynomial $P$ and the natural normal conditional expectation 
$\mathcal{E}_G :\Big(\prod_{\U}L_1(N_m)\Big)^* \rightarrow \pi_G(VN(\mathbb{F}_\infty))$,
which is the adjoint of the inclusion $\pi_G(VN(\mathbb{F}_\infty))\hookrightarrow \Big(\prod_{\U}L_1(N_m)\Big)^*$.

Combining the above discussions we have the following representations of vector valued noncommuatative $L_p$ spaces 
with respect to $$\N := VN(\mathbb{F}_\infty)\,\, \text{and}\,\, N := \{\WT_i : i\geq 1\}'',$$ which are both satisfying $QWEP$. 
These $\N$, $N$, $\N_m$ and $N_m$ will be fixed from now on.

\begin{prop}\label{prop-repVecLpFreeGpVN}
Let $E \subseteq C^*(\mathbb{F}_{\infty})$, $1\leq p < \infty$ and $\U$ be a free ultrafilter in $\n$. 
Then, for any $n \in \mathbb{N}$ and $(x_i)^n_{i=1} \subseteq E$, we have
$$\norm{\sum^n_{i=1} \lambda(g_i) \otimes x_i}_{L_p(\N ; E)} = 
\lim_{\U} \norm{\sum^n_{i=1}U^m_i \otimes x_i}_{L_p(\N_m ;E)}$$ and 
$$\norm{\sum^n_{i=1}\widetilde{W_i} \otimes x_i}_{L_p(N ; E)}
 = \lim_{\U} \norm{\sum^n_{i=1}G^m_i \otimes x_i}_{L_p(N_m ;E)}.$$
\end{prop}

It would be convenient for later use to fix $\alpha_U = (\pi_U, \mathcal{E}_U)$ 
and $\alpha_G = (\pi_G, \mathcal{E}_G)$ for $\N$ and $N$, respectively, with the same free ultrafilter $\U$ on $\n$ 
and use the notation $L_p(\N ; E)$ and $L_p(N ; E)$ for $L_p(\N, \alpha_U ; E)$ and $L_p(N, \alpha_G ; E)$.

\section{Free probabilistic tools}\label{sec-tools}

From now on, $\M$ and $E$ implies a von Neumann algebra satisfying $QWEP$ 
and a locally-$C^*(\mathbb{F}_{\infty})$ operator space, respectively. 
The reasons for this restriction comes from the discussions in the subsection \ref{subsec-LpQWEP}.
Moreover, we fix $\alpha = (\pi, \mathcal{E})$ for $\M$.

In this section we collect operator space analogues of several probabilistic tools in Banach space theory.
We start with the ``contraction principle for $*$-free Haar Unitaries".

\begin{thm}\label{thm-contraction-free}
Let $(\phi_i)_{i\geq 1}$ be a sequence in $L_p(\M)$ $(1 \leq p < \infty)$ 
with $$\norm{\sup_i \abs{\phi_i}}_p < \infty.$$ 
Then we have for any $n \in \mathbb{N}$ and $(x_i)^n_{i=1} \subseteq E$
$$\norm{\sum^n_{i=1} \phi_i \otimes \lambda(g_i) \otimes x_i}_{L_p(\M \overline{\otimes} \N ; E)} 
\leq 2d_f(E)\norm{\sup_i \abs{\phi_i}}_p\norm{\sum^n_{i=1}\lambda(g_i) \otimes x_i}_{L_p(\N ; E)}.$$ 
\end{thm}
\begin{proof}
We assume that $E \subseteq C^*(\mathbb{F}_{\infty})$. By a usual density argument we can also assume that $\M$ is $\sigma$-finite. 
Then there is a normal faithful state $\varphi$ on $\M$. 
Now we suppose that $$\norm{\sup_i \abs{\phi_i}}_p = 1$$ and fix $n \in \n$.
For any given $\epsilon > 0 $ we can find $a, b \in L_{2p}(\M)$ 
and contractions $(x_i)_{i\geq 1} \subseteq \M$ such that $$\phi_i = a x_i b\,\, \text{for all} \,\, i\geq 1\,\, 
\text{and}\,\, \norm{a}_{L_{2p}(\M)} = 1,\,\, \norm{b}_{L_{2p}(\M)} \leq 1+\epsilon.$$
By a standard argument we can find $D (\geq 0) \in L_1(\M)$ and contractions $(u_i)_{i\geq 1} \subseteq \M$ 
such that $$\phi_i = D^{\frac{1}{2p}} u_i D^{\frac{1}{2p}}\,\, \text{for all} \,\, i\geq 1\,\, \text{and} 
\,\, \norm{D}_{L_1(\M)} \leq 4(1+\epsilon'),$$ where $\epsilon' \rightarrow 0$ as $\epsilon \rightarrow 0.$
Indeed, if we set $D = (\abs{a}^p + \abs{b}^p)^2$ we have $\abs{a}, \abs{b} \leq D^{\frac{1}{2p}}$ 
so that there are contractions $V,W \in \M$ such that $$a = D^{\frac{1}{2p}}V \,\, \text{and} \,\, b = WD^{\frac{1}{2p}}.$$ 
Then we have $\phi_i = a x_i b = D^{\frac{1}{2p}}V x_i WD^{\frac{1}{2p}} = D^{\frac{1}{2p}} u_i D^{\frac{1}{2p}}$ 
with $u_i = V x_i W$, and $$\norm{D}_{L_1(\M)} \leq \norm{a}^{2p}_{L_{2p}(\M)} + \norm{b}^{2p}_{L_{2p}(\M)} 
+ 2\norm{\abs{a}^p\abs{b}^p}_{L_1(\M)} \leq 1 + (1+\epsilon)^{2p} + 2 (1+\epsilon)^2.$$
By Russo-Dye theorem and the usual convexity argument it is enough to show the case that $u_i$'s are all unitaries.

With all these assumption for any $(x_i)^n_{i=1} \subseteq E$ we have 
$$\sum^n_{i=1}\phi_i \otimes \lambda(g_i) \otimes x_i = \sum^n_{i=1} D^{\frac{1}{2p}} u_i D^{\frac{1}{2p}} \otimes \lambda(g_i) \otimes x_i 
= \Phi\circ \pi \otimes I_E \Big(\sum^n_{i=1}\lambda(g_i) \otimes x_i\Big),$$ where 
$$\Phi : L_p(\M\overline{\otimes}\N) \rightarrow L_p(\M\overline{\otimes}\N),\,\, 
z \mapsto d^{\frac{1}{2p}} z d^{\frac{1}{2p}}$$ with $d = D  \otimes I_{\N}$ and 
$$\pi :  L_p(\N) \rightarrow  L_p(\M\overline{\otimes}\N)$$ is the map induced from the $*$-isomorphism 
$$\pi_{\infty} : \N \rightarrow \M\overline{\otimes}\N,\,\, \text{with}\,\, \pi_{\infty}(\lambda(g_i)) = u_i \otimes \lambda(g_i).$$
Note that it is easy to check that $(\lambda(g_i))_{i\geq 1}$ and $(u_i \otimes \lambda(g_i))_{i\geq 1}$ 
have the same $*$-distribution on $W^*$-probability spaces with faithful states $(\N, \tau)$ 
and $(\M\overline{\otimes}\N, \varphi \otimes \tau)$, respectively. Thus, $\pi_{\infty}$ is a $*$-isomorphism by Lemma 1 of \cite{HiPe00}.

Now we have $\Phi$ and $\pi$ are completely positive, and consequently so is the composition $\Phi\circ \pi$ with
\begin{align*}
\begin{split}
\norm{\Phi\circ \pi} & = \norm{\Phi\circ \pi(I)} = 
\norm{d^{\frac{1}{2}}}_{L_2(\M\overline{\otimes}\N)} = \norm{D^{\frac{1}{2}}}_{L_2(\M)}\\ 
& = \norm{D}^{\frac{1}{2}}_{L_1(\M)} \leq 2\sqrt{1+\epsilon'}.
\end{split}
\end{align*}
Then, by the Proposition \ref{prop-CPmap} we have 
\begin{align*}
\begin{split}
\norm{\sum^n_{i=1} \phi_i \otimes \lambda(g_i) \otimes x_i}_{L_p(\M \overline{\otimes} \N ; E)} 
& \leq \norm{\Phi\circ \pi} \norm{\sum^n_{i=1}\lambda(g_i) \otimes x_i}_{L_p(\N ; E)}\\
& \leq 2\sqrt{1+\epsilon'} \norm{\sum^n_{i=1}\lambda(g_i) \otimes x_i}_{L_p(\N ; E)}.
\end{split}
\end{align*}

\end{proof}

Now we consider a partial relationship between the average with respect to $*$-free Haar unitaries 
and the average with respect to $*$-free circular elements. 
In the proof we will use the random matrix model approach in the subsection \ref{subsec-RandomMatrix}.

\begin{thm}\label{thm-free-Haar-circular}
Let $1\leq p < \infty$. Then there is a constant $C_G > 0$ such that 
we have for any $n \in \mathbb{N}$ and $(x_i)^n_{i=1} \subseteq E$
$$\norm{\sum^n_{i=1} \lambda(g_i) \otimes x_i}_{L_p(\N ; E)} 
\leq C_G\norm{\sum^n_{i=1}\widetilde{W}_i \otimes x_i}_{L_p(N ; E)}.$$ 
\end{thm}
\begin{proof}
First, we assume that $E \subseteq C^*(\mathbb{F}_{\infty})$, and fix $n \in \mathbb{N}$ and $(x_i)^n_{i=1} \subseteq E$. 
Let $(U^m_i)^n_{i=1}$ and $(G^m_i)^n_{i=1}$ be the independent familys of the standard unitary random matrices 
and the standard gaussian random matrices as in the subsection \ref{subsec-RandomMatrix}. 
Note that $(U^m_i\abs{G^m_i})^n_{i=1}$ and $(G^m_i)^n_{i=1}$ have the same $*$-distribution and 
$$\E(U^m_i\abs{G^m_i}) = U^m_i\E\abs{G^m_i},$$ where $\E$ implies the expectation with respect to $(\Om', P')$. 
Moreover, it is well known that (Remark 1.7. in p.80 of \cite{MP81}) $$\E\abs{G^m_i} = \delta(m) I_m,$$ 
where $I_m$ is the identity matrix on $M_m$ and $\delta(m) \geq \delta$ for some constant $\delta > 0$. 
Thus, for $\N_m = L_\infty(\Om; M_m)$  and $N_m = L_\infty(\Om'; M_m)$, we have 
\begin{align*}
\begin{split}
\norm{\sum^n_{i=1}U^m_i \otimes x_i}_{L_p(\N_m ;E)} & = 
\delta(m)^{-1} \norm{\sum^n_{i=1}U^m_i\E\abs{G^m_i} \otimes x_i}_{L_p(\N_m ;E)}\\
& = \delta(m)^{-1} \norm{\sum^n_{i=1}\E(U^m_i\abs{G^m_i}) \otimes x_i}_{L_p(\N_m ;E)}\\
& \leq \delta^{-1} \norm{\sum^n_{i=1}U^m_i\abs{G^m_i} \otimes x_i}_{L_p(\Om \times \Om' ; M_m(E))}\\
& = \delta^{-1} \norm{\sum^n_{i=1}G^m_i \otimes x_i}_{L_p(N_m ;E)}
\end{split}
\end{align*}
Let $\U$ be the free ultrafilter on $\n$ chosen for $\alpha_U = (\pi_U, \mathcal{E}_U)$ 
and $\alpha_G = (\pi_G, \mathcal{E}_G)$, then we have
\begin{align*}
\begin{split}
\norm{\sum^n_{i=1} \lambda(g_i) \otimes x_i}_{L_p(\N ; E)} & = 
\lim_{m, \U} \norm{\sum^n_{i=1}U^m_i \otimes x_i}_{L_p(\N_m ;E)}\\
& \leq \delta^{-1}\lim_{m, \U} \norm{\sum^n_{i=1}G^m_i \otimes x_i}_{L_p(N_m ;E)}\\
& = \delta^{-1} \norm{\sum^n_{i=1}\widetilde{W_i} \otimes x_i}_{L_p(N ; E)}
\end{split}
\end{align*}

\end{proof}

Another important ingredient is the concept of $K$-convexity with respect to circular elements.

\begin{defn}
Let $E$ be a locally-$C^*(\mathbb{F}_{\infty})$ operator space. 
Then we define the ``$\WT K$-convexity" constant $\WT K(E)$ by $\WT K(E) = \sup_n \WT K_n(E),$
where $$\WT K_n(E) = \norm{P_n \otimes I_E : L_2(N;E) \rightarrow L_2(N;E)}$$ and
$$P_n : L_2(N) \rightarrow L_2(N), \,\, a \mapsto \sum^n_{i=1}\left\langle a, \widetilde{W_i} \right\rangle \widetilde{W_i}.$$
\end{defn}

\begin{thm}\label{thm-K-convex-circular}
Suppose $S_2(E)$ is K-convex as a Banach space.
Then $E$ is ``$\WT K$-convex" with $$\WT K(E) \leq K(S_2(E)),$$ 
where $K(X)$ is the $K$-convexity constant of a Banach space $X$.
\end{thm}
\begin{proof} 
Using the isometric embedding $$L_2(N;E) \hookrightarrow \prod_{m, \U}L_2(N_m ;E),$$ 
where $N_m = L_{\infty}(\Om' ;M_m)$, we will take a detour down to $m$-th random matrix level. 

Let $$\G^n_{m} = \text{span}\{ g^k_{ij} : 1\leq k \leq n , \,\, 1\leq i,j \leq m \}\subseteq L_2(\Om')$$
and $$\widetilde{\G}_{m} = \text{span}\{G^m_k : 1\leq k \leq n\} \subseteq L_2(N_m).$$
Let $Q_m$ be the orthogonal projection from $L_2(\Om')$ onto $\G^n_{m}$
and $$R_m : \G^n_{m}(L_2(M_m)) \rightarrow \widetilde{\G}_{m}$$ be the linear map defined by
$$R_m \Big(\sum^n_{k=1} \sum^m_{i,j =1} g^k_{ij} \sum^m_{r,s =1}e_{rs} \otimes x^{ijk}_{rs} \Big) =
\sum^n_{k=1} \Big(\sum^m_{i,j =1} g^k_{ij} e_{ij} \Big) \otimes \text{Av}(x^{ijk}_{rs})$$ for 
any $x^{ijk}_{rs} \in \mathbb{C}\,\, (1 \leq i,j,r,s \leq m, 1 \leq k \leq n)$, 
where $$\text{Av}(x^{ijk}_{rs}) = \frac{1}{\abs{\Sigma_m}^2}\sum_{\sigma, \rho \in \Sigma_m} y^{k}_{\sigma(i)\rho(j)}, 
\,\, y^{k}_{ij} = x^{ijk}_{ij}$$ and $\Sigma_m$ is the set of all permutations of $\{1, \cdots, m\}$.
Note that $\text{Av}(x^{ijk}_{rs})$ depends only on $k$. 

Then we can recover $P_n$ by taking ultraproduct of $R_m\circ (Q_m \otimes I_{M_m})$.
More precisely, for any $f \in L_2(N)\otimes E$ we can associate $(f_m)_{m,\U} \in L_2(N_m)$
obtained from the isometric embedding $L_2(N) \hookrightarrow \prod_{m,\U}L_2(N_m)$. 
Then we have $$P_n(f) = \lim_{m,\U}R_m\circ (Q_m\otimes I_{M_m})(f_m).$$
Thus, it is enough to show that $$Q_m \otimes I_{M_m} \otimes I_E : L_2(\Om')(M_m(E))\rightarrow \G^n_{m}(M_m(E))$$ 
and $$R_m \otimes I_E : \G^n_{m}(L_2(M_m ; E)) \rightarrow \widetilde{\G}_{m}(E)$$ are both bounded. 
The first one is obtained by the assumption of $K$-convexity of $S_2(E)$ so that we have 
$$\norm{Q_m \otimes I_{M_m} \otimes I_E} \leq K(S^m_2(E)) \leq K(S_2(E)).$$

Now we consider the map $R_m \otimes I_E$. Let $(\varepsilon_i)_{i\geq1}$ is the Rademacher sequences on $[0,1]$. 
Then since $(g^k_{ij})$ and $(\epsilon^k_{ij}g^k_{ij})$ 
have the same distribution for any $\epsilon^k_{ij} \in \{\pm 1\}$ we have 
\begin{align*}
\begin{split}
\lefteqn{\norm{\sum^n_{k=1} \sum^m_{i,j =1} g^k_{ij} \sum^m_{r,s =1}e_{rs} \otimes x^{ijk}_{rs}}_{L_2(N_m; E)}}\\
& = \int^1_0\norm{\sum^n_{k=1} \sum^m_{i,j,r,s =1} \varepsilon_i(t) \varepsilon_r(t)
g^k_{ij} e_{rs} \otimes x^{ijk}_{rs}}_{L_2(N_m; E)}dt\\
& \geq \norm{\sum^n_{k=1} \sum^m_{i,j,r,s =1} \Big[\int^1_0 \varepsilon_i(t) \varepsilon_r(t)dt\Big]
g^k_{ij} e_{rs} \otimes x^{ijk}_{rs}}_{L_2(N_m; E)}\\
& = \norm{\sum^n_{k=1} \sum^m_{i,j,s =1} g^k_{ij} e_{is} \otimes x^{ijk}_{is}}_{L_2(N_m; E)}.
\end{split}
\end{align*}
By repeating the same procedure we get 
$$\norm{\sum^n_{k=1} \sum^m_{i,j =1} g^k_{ij} \sum^m_{r,s =1}e_{rs} \otimes x^{ijk}_{rs}}_{L_2(N_m; E)}
\geq \norm{\sum^n_{k=1} \sum^m_{i,j =1} g^k_{ij} e_{ij} \otimes y^{k}_{ij}}_{L_2(N_m; E)}$$ 
for $y^{k}_{ij} = x^{ijk}_{ij}$. 

We can proceed further using permutations. 
Note that $(g^k_{ij})$ and $(g^k_{\sigma^{-1}(i)j})$ 
have the same distribution for any permutation $\sigma \in \Sigma_m$, and also 
interchanging $i$-th row into $\sigma^{-1}(i)$-th row does not affect the norm. Thus we have
\begin{align*}
\begin{split}
\lefteqn{\norm{\sum^n_{k=1} \sum^m_{i,j =1} g^k_{ij} e_{ij} \otimes 
\frac{1}{\abs{\Sigma_m}} \sum_{\sigma \in \Sigma_m}y^k_{\sigma(i)j}}_{L_2(N_m; E)}}\\ & \leq
\frac{1}{\abs{\Sigma_m}} \sum_{\sigma \in \Sigma_m}\norm{\sum^n_{k=1} \sum^m_{i,j =1} 
g^k_{ij} e_{ij} \otimes y^k_{\sigma(i)j}}_{L_2(N_m; E)}\\ & =
\frac{1}{\abs{\Sigma_m}} \sum_{\sigma \in \Sigma_m}\norm{\sum^n_{k=1} \sum^m_{i,j =1} 
g^k_{\sigma^{-1}(i)j} e_{\sigma^{-1}(i)j} \otimes y^k_{ij}}_{L_2(N_m; E)}\\ & =
\frac{1}{\abs{\Sigma_m}} \sum_{\sigma \in \Sigma_m}\norm{\sum^n_{k=1} \sum^m_{i,j =1} 
g^k_{ij} e_{\sigma^{-1}(i)j} \otimes y^k_{ij}}_{L_2(N_m; E)}\\ & =
\frac{1}{\abs{\Sigma_m}} \sum_{\sigma \in \Sigma_m}\norm{\sum^n_{k=1} \sum^m_{i,j =1} 
g^k_{ij} e_{ij} \otimes y^k_{ij}}_{L_2(N_m; E)}\\ & = \norm{\sum^n_{k=1} \sum^m_{i,j =1} 
g^k_{ij} e_{ij} \otimes y^k_{ij}}_{L_2(N_m; E)}.
\end{split}
\end{align*}
By repeating the same procedure we get $$ \norm{\sum^n_{k=1} \sum^m_{i,j =1} g^k_{ij} e_{ij} \otimes y^{k}_{ij}}_{L_2(N_m; E)}
\geq \norm{\sum^n_{k=1} \sum^m_{i,j =1} g^k_{ij} e_{ij} \otimes \text{Av}(x^{ijk}_{rs})}_{L_2(N_m; E)},$$ 
and combining these results we get that $R_m \otimes I_E$ is a contraction.

\end{proof}

In Banach space theory, one of the useful characterizations of Hilbert spaces is being 
type 2 and cotype 2 simultaneously. We will need similar characterization of operator Hilbert spaces.

First, we define a noncommutative analogue of type 2 and cotype 2 using 
completely 2-summing norm and a variant of ``$\ell$-norm". Note that this is a variant of $S_2$-type and $S_2$-cotype defined in \cite{L-typecotypeOS}. 
Recall that for a linear map $u : E\rightarrow F$ between operator spaces we define 
$$\pi^o_2(u) = \sup\Big\{ \norm{(ux_{ij})}_{S^n_2(F)} : 
\norm{\sum^n_{i,j=1}e_{ij}\otimes x_{ij}}_{S^n_2 \otimes_{\min} E} \leq 1,\,\, n\in \n\Big\}.$$

\begin{defn}
Let $u : \ell^n_2 \rightarrow E$, $n\in \n$. 
Then we define ``$\ell_{\WT}$-norm" of $u$ by 
$$\ell_{\WT}(u) := \norm{\sum^n_{i=1}\widetilde{W}_i \otimes ue_i}_{L_2(N ;E)}.$$
\end{defn}

\begin{defn}
We say that $E$ has $\WT$-type 2 if there is a constant $C>0$ such that 
$$\ell_{\WT}(u) \leq C \pi^o_2(u^*)$$ 
for any $n\in \n$ and $u : OH_n \rightarrow E$.

We say that $E$ has $\WT$-cotype 2 if there is a constant $C'>0$ such that 
$$\pi^o_2(u) \leq C'\ell_{\WT}(u)$$ 
for any $n\in \n$ and $u : OH_n \rightarrow E$. 
We denote ${\WT T}^{oh}_2(E)$ and ${\WT C}^{oh}_2(E)$ for the infimum of such $C$ and $C'$, respectively.
\end{defn}

\begin{prop}\label{prop-OH-TypeCotypeCharacterization}
$E$ is completely isomorphic to an operator Hilbert space if and only if $E$ has $\WT$-type 2 and $\WT$-cotype 2. 
In this case we have $$d_{cb}(E, OH(I)) \leq {\WT T}^{oh}_2(E){\WT C}^{oh}_2(E)$$ for some index set $I$.
\end{prop}
\begin{proof}
By applying trace duality to (iii) of Theorem 6.5 \cite{P98} 
it is enough to show that $$\pi^o_2(u) \leq {\WT T}^{oh}_2(E){\WT C}^{oh}_2(E)\pi^o_2(u^*)$$ for any $n\in \n$ 
and $u : OH_n \rightarrow E$. Indeed, we have 
$$\pi^o_2(u) \leq {\WT C}^{oh}_2(E)\ell_{\WT}(u) \leq {\WT T}^{oh}_2(E){\WT C}^{oh}_2(E)\pi^o_2(u^*).$$
\end{proof}

\begin{prop}\label{prop-ell-norm}
Let $u : \ell^n_2 \rightarrow E$, 
$v: E \rightarrow \ell^n_2$ and $A : \ell^n_2 \rightarrow \ell^n_2$ for $n\in \n$. 
Then we have
\begin{itemize}
\item[(1)]
$\ell_{\WT}(uA) \leq \ell_{\WT}(u)\norm{A}.$
\item[(2)]
$\ell_{\WT}(v^*) \leq d_f(E)\WT K(E)\ell^*_{\WT}(v),$
\end{itemize}
where $\ell^*_{\WT}(\cdot)$ refers to the trace dual norm of $\ell_{\WT}(\cdot)$.
\end{prop}
\begin{proof}
(1) We consider the random matrix model again. By the usual convexity argument it is enough to show that for any $m\in \n$
$$\norm{\sum^n_{i=1}G^m_i \otimes uAe_i}_{L_2(N_m;E)}  = \norm{\sum^n_{i=1}G^m_i \otimes ue_i}_{L_2(N_m;E)}$$ 
whenever $A = (a_{ij})$ is unitary. Indeed, we have 
$$\sum^n_{i=1}G^m_i \otimes uAe_i = \sum^n_{i=1} G^m_i \otimes u\Big(\sum^n_{j=1}a_{ji}e_j \Big) 
= \sum^n_{j=1}\Big(\sum^n_{i=1}a_{ji}G^m_i \Big) \otimes ue_j$$
and since $(G^m_j)^n_{j=1}$ and $(\sum^n_{i=1}a_{ji}G^m_i)^n_{j=1}$ have the same $*$-distribution we are done.

(2) Consider $\sum^n_{i=1}\WT_i \otimes v^*e_i \in L_2(N;E^*)$. 
Since $L_2(N ; E^*) \hookrightarrow L_2(N ; E)^*$ completely isomorphically (Proposition \ref{prop-duality})
for any $\epsilon > 0$ we can choose $f \in L_2(N ; E)$ with $$\norm{f}_{L_2(N ; E)} =1\,\, \text{and}
\,\, \abs{\left\langle \sum^n_{i=1}\WT_i \otimes v^*e_i, f \right\rangle} \geq d_f(E)^{-1}(1- \epsilon) \ell_{\WT}(v^*).$$
Then for $w : \ell^n_2 \rightarrow E,\,\, e_i \rightarrow \left\langle \WT_i, f \right\rangle$ we have
\begin{align*}
\begin{split}
d_f(E)^{-1}(1- \epsilon) \ell_{\WT}(v^*)& \leq \abs{\left\langle \sum^n_{i=1}\WT_i \otimes v^*e_i, f \right\rangle} 
= \abs{\text{tr}(vw)}\\ & \leq \ell_{\WT}(w)\ell^*_{\WT}(v)\leq \WT K(E) \norm{f}_{L_2(N ; E)} \ell^*_{\WT}(w)\\
& = \WT K(E) \ell^*_{\WT}(w).
\end{split}
\end{align*}
Since $\epsilon>0$ is arbitrary we get the desired result.

\end{proof}

\begin{rem}{\rm
Let $(\widetilde{W}_{ij})_{i,j\geq 1}$ be a re-indexing of $(\widetilde{W}_{i})_{i\geq 1}$. 
Then by (1) of Proposition \ref{prop-ell-norm} we can easily check that $E$ has $\WT$-type 2 
if and only if there is a constant $C>0$ such that 
$$\norm{\sum^n_{i,j=1}\widetilde{W}_{ij}\otimes x_{ij}}_{L_2(N;E)} \leq C \norm{(x_{ij})}_{S^n_2(E)}$$ 
for any $n\in \n$ and $x_{ij} \in E$. 

Similarly, $E$ has $\WT$-cotype 2 if and only if there is a constant $C'>0$ such that 
$$\norm{(x_{ij})}_{S^n_2(E)} \leq C'\norm{\sum^n_{i,j=1}\widetilde{W}_{ij}\otimes x_{ij}}_{L_2(N;E)}$$ 
for any $n\in \n$ and $x_{ij} \in E$.
}\end{rem}

As in the Banach space case we have a duality of $\W$-type 2 for ``$\WT K$-convex" operator space. 

\begin{prop}\label{prop-type-duality}
${\WT T}^{oh}_2(E) \leq d_f(E){\WT K}(E){\WT C}^{oh}_2(E^*).$
\end{prop}
\begin{proof}
Note that $\pi^o_2$ is self-dual in the sense of trace duality. (See \cite{L-typecotypeOS} for example) 
By applying (2) of Proposition \ref{prop-ell-norm} and trace duality we have
\begin{align*}
\begin{split}
{\WT T}^{oh}_2(E) & = \sup\Big\{ \frac{\pi^o_2(v^*)}{\ell^*_{\WT}(v)}: n\in \n\,\, \text{and}\,\, v : E \rightarrow OH_n \Big\}\\
& \leq d_f(E){\WT K}(E) \cdot \sup\Big\{ \frac{\pi^o_2(v^*)}{\ell_{\WT}(v^*)}: n\in \n\,\, \text{and}\,\, v : E \rightarrow OH_n \Big\}\\
& = d_f(E){\WT K}(E){\WT C}^{oh}_2(E^*).
\end{split}
\end{align*}

\end{proof}

\section{Unconditionality with respect to complete orthonormal systems in noncommutative $L_2$ spaces}

We start this section with an equivalent formulation of completely 2-summing norm.
The following proposition tells us that we can replace $(e_{ij})$ into complete orthonormal systems 
in some noncommutative $L_2$-space.
\begin{prop}\label{prop-replacing-2summing}
Let $F$ be a locally-$C^*(\mathbb{F}_{\infty})$ operator space, $\M$ be a von Neumann algebra not subhomogeneous, 
$\Phi=(\phi_i)_{i\geq 1}$ be a complete orthonormal system in $L_2(\M)$ and $u : E\rightarrow F$. 
Then $u$ is completely 2-summing if and only if there is a constant $C>0$ such that
\begin{eqnarray}\label{eqn-com-2summing}
\norm{\sum^n_{i,j=1}\phi_{ij}\otimes ux_{ij}}_{L_2(\M;F)}
\leq C\norm{\sum^n_{i,j=1}e_{ij}\otimes x_{ij}}_{S^n_2 \otimes_{\min} E}
\end{eqnarray}
for any $n\in \n$ and $(x_{ij}) \in M_n(E)$, where $(\phi_{ij})_{i,j \geq 1}$ be a re-indexing of $(\phi_i)_{i\geq 1}$.
Moreover, $$\frac{1}{d_f(F)}\cdot \inf C \leq \pi^o_2(u) \leq d_f(F)\cdot \inf C,$$ 
where the infimum runs over the $C$'s satisfying (\ref{eqn-com-2summing}).
\end{prop}
\begin{proof}
First, we assume that $F \subseteq C^*(\mathbb{F}_{\infty})$.

Now suppose (\ref{eqn-com-2summing}) holds, and fix $n \in \n$ and $(x_{ij}) \in M_n(E)$. 
Since $\M$ is not subhomogeneous for any $\epsilon>0$ there are completely positive maps
$$\rho : M_n \rightarrow \M \,\, \text{and}\,\, \sigma : \M \rightarrow M_n\,\, \text{such that}\,\, 
\norm{\sigma \rho - I_{M_n}}_{cb} \leq \epsilon$$ by Lemma 2.7 of \cite{S83} and \cite{HT83}. 
Thus by Proposition \ref{prop-CPmap} we can replace matrix units $(e_{ij})^n_{i, j=1} \in M_n$ into 
for some $(f_{ij})^n_{i, j=1} \in \M$ allowing constant $1+ \epsilon$.
Since $\Phi$ is complete in $L_2(\M)$ we can choose $(y_{kl})^m_{k,l=1} \subseteq E$ 
such that $$\norm{\sum^n_{i,j=1}f_{ij}\otimes x_{ij} - \sum^m_{k,l=1}\phi_{kl}\otimes y_{kl}}_{L_2(\M;F)} \leq \epsilon.$$
Then we have 
\begin{align*}
\begin{split}
\norm{(ux_{ij})}_{S^n_2(F)} & \leq (1+\epsilon) \norm{\sum^n_{i,j=1}f_{ij}\otimes u x_{ij}}_{L_2(\M ;E)}\\ 
& \leq (1+\epsilon)\norm{u}_{cb}\norm{\sum^n_{i,j=1}f_{ij}\otimes x_{ij} - \sum^m_{k,l=1}\phi_{kl}\otimes y_{kl}}_{L_2(\M;F)}
\\ & \;\;\;\; + (1+\epsilon)\norm{\sum^m_{k,l=1}\phi_{kl}\otimes uy_{kl}}_{L_2(\M;F)}\\
& \leq \epsilon(1+\epsilon) \norm{u}_{cb} + C(1+\epsilon)\norm{(y_{kl})}_{S^m_2 \otimes_{\min}E}
\end{split}
\end{align*}
and
\begin{align*}
\begin{split}
\norm{(y_{kl})}_{S^m_2\otimes_{\min}E} & = 
\norm{E^* \rightarrow S^m_2,\,\, e^*\mapsto (\left\langle e^*, y_{kl} \right\rangle)}_{cb}\\
& = \norm{S_2(E^*) \rightarrow S_2(S^m_2),\,\, e^*_{rs}\mapsto 
(\left\langle e^*_{rs}, y_{kl} \right\rangle)_{1\leq r,s,\, 1\leq k,l\leq m}}\\
& = \sup_{\norm{(e^*_{rs})}_{S_2(E^*)} \leq 1} \norm{(\left\langle e^*_{rs}, y_{kl} \right\rangle)}_{S_2(S^m_2)}\\
& = \sup_{\norm{(e^*_{rs})}_{S_2(E^*)} \leq 1} 
\norm{I_{L_2(\M)} \otimes \Psi_{(e^*_{rs})}\Big(\sum^m_{k,l=1}\phi_{kl}\otimes y_{kl}\Big)}_{L_2(\M ;S_2)},
\end{split}
\end{align*}
where $$\Psi_{(e^*_{rs})} : E \rightarrow S_2, \,\, 
z \mapsto \Big(\left\langle e^*_{rs}, z \right\rangle \Big)_{r, s}.$$
Note that we have by (2) of Proposition \ref{prop-LpQWEP-Basic} and Lemma 5.14. of \cite{P98} that 
$$\norm{I_{L_2(\M)} \otimes \Psi_{(e^*_{rs})}}_{cb} \leq \norm{\Psi_{(e^*_{rs})}}_{cb} 
\leq \pi^o_2(\Psi_{(e^*_{rs})}) \leq \norm{(e^*_{rs})}_{S_2(E^*)}.$$
Since $(f_{ij})$ is orthonormal we have
\begin{align*}
\begin{split}
\lefteqn{\norm{(y_{kl})}_{S^m_2\otimes_{\min}E}}\\ & \leq \sup_{\norm{(e^*_{rs})}_{S_2(E^*)} \leq 1}
\norm{I_{L_2(\M)} \otimes \Psi_{(e^*_{rs})}\Big(\sum^n_{i,j=1}f_{ij}\otimes x_{ij} 
- \sum^m_{k,l=1}\phi_{kl}\otimes y_{kl}\Big)}_{L_2(\M ;S_2)} \\ & \,\,\,\,\, + 
\sup_{\norm{(e^*_{rs})}_{S_2(E^*)} \leq 1} 
\norm{I_{L_2(\M)} \otimes \Psi_{(e^*_{rs})}\Big(\sum^n_{i,j=1}f_{ij}\otimes x_{ij}\Big)}_{L_2(\M ;S_2)}\\
& \leq \norm{\sum^n_{i,j=1}f_{ij}\otimes x_{ij} - \sum^m_{k,l=1}\phi_{kl}\otimes y_{kl}}_{L_2(\M ;E)}
+ \sup_{\norm{(e^*_{rs})}_{S_2(E^*)} \leq 1} \norm{(\left\langle e^*_{rs}, x_{ij} \right\rangle)}_{S_2(S^n_2)}\\
& \leq \epsilon + \norm{E^* \rightarrow S^n_2,\,\, e^*\mapsto (\left\langle e^*, x_{ij} \right\rangle)}_{cb} 
= \epsilon + \norm{(x_{ij})}_{S^n_2\otimes_{\min}E}.
\end{split}
\end{align*}
By combining all these results we get
$$\norm{(ux_{ij})}_{S^n_2(F)} \leq \epsilon(1+\epsilon) \norm{u}_{cb} + C\epsilon(1+\epsilon) + C(1+\epsilon)\norm{(x_{ij})}_{S^n_2\otimes_{\min}E},$$
and by letting $\epsilon \rightarrow 0$ we get
$\norm{(ux_{ij})}_{S^n_2(F)} \leq C\norm{(x_{ij})}_{S^n_2\otimes_{\min}E}$ and consequently
$$\pi^o_2(u) \leq C.$$

For the converse we consider a completely 2-summing map $u : E\rightarrow F$.
By the factorization theorem (Proposition 6.1. of \cite{P98}) we have 
$$A : E \rightarrow OH(I)\,\, \text{and}\,\, B : OH(I) \rightarrow F$$ for some index set $I$
such that $u = BA$ with $\pi^o_2(A)\leq 1$ and $\norm{B}_{cb}\leq \pi^o_2(u)$.

Then by (2) of Proposition \ref{prop-LpQWEP-Basic} we have 
\begin{align*}
\begin{split}
\norm{\sum^n_{i,j=1}\phi_{ij}\otimes ux_{ij}}_{L_2(\M;E)} & \leq 
\norm{B}_{cb}\norm{\sum^n_{i,j=1}\phi_{ij}\otimes Ax_{ij}}_{L_2(\M;E)}\\
& \leq \pi^o_2(u)\norm{(Ax_{ij})}_{S^n_2(OH(I))}\\
& \leq \pi^o_2(u)\pi^o_2(A)\norm{(x_{ij})}_{S^n_2 \otimes_{\min}E}\\
& \leq \pi^o_2(u)\norm{(x_{ij})}_{S^n_2 \otimes_{\min}E}
\end{split}
\end{align*}
for any $n\in \n$ and $(x_{ij}) \subseteq E$.

\end{proof}

Now we define the unconditionality with respect to orthonormal systems in noncommutative $L_2$ spaces.

\begin{defn} 
Let $\Phi=(\phi_i)_{i\geq 1}$ be an orthonormal system in $L_2(\M)$. We say that $E$ is $\Phi$-unconditional if 
there is a constant $C>0$ such that
$$\frac{1}{C} \norm{\sum^n_{i=1} \phi_i \otimes x_i}_{L_2(\M;E)}\leq 
\norm{\sum^n_{i=1} \phi_i \otimes u_i \otimes x_i}_{L_2(\M \overline{\otimes}M_m ;E)} 
\leq C \norm{\sum^n_{i=1} \phi_i \otimes x_i}_{L_2(\M;E)}$$ for any $n,m \in \n$,
$(x_i)^n_{i=1}\subseteq E$ and any unitaries $u_i \in M_m$. We denote $\Phi_{unc}(E)$ by the infimum of such $C$.
\end{defn}

Before we prove our main result we observe that unconditionality implies $\WT$-cotype 2.
\begin{prop}\label{prop-UncImpliesW-cotype2}
Let $\M$ be a von Neumann algebra not subhomogeneous, and $\Phi=(\phi_i)_{i\geq 1}$ be
a complete orthonormal system in $L_2(\M)$ with $$\norm{\sup_i \abs{\phi_i}}_2 < \infty.$$ 
If $E$ is $\Phi$-unconditional, then $E$ has $\WT$-cotype 2 with
$${\WT C}^{oh}_2(E) \leq 2C_G d_f(E)^2\norm{\sup_i \abs{\phi_i}}_2\Phi_{unc}(E),$$ 
where $C_G$ is the constant in Theorem \ref{thm-contraction-free}. 
\end{prop}
\begin{proof}
Let's fix $k, n\in \n$, $(x_{ij})^n_{i,j=1} \subseteq OH_k$ and $u : OH_k \rightarrow E$ and let 
$(\phi_{ij})_{i,j \geq 1}$ be a re-indexing of $(\phi_i)_{i\geq 1}$.

It is clear from the definition that 
$$\norm{\sum^n_{i,j=1} \phi_{ij} \otimes ux_{ij}}_{L_2(\M;E)}\leq \Phi_{unc}(E)
\norm{\sum^n_{i=1} \phi_{ij} \otimes U^m_{ij} \otimes ux_{ij}}_{L_2(\M \overline{\otimes} N_m ;E)}$$
for independent family of standard random unitaries $U^m_{ij}$, $m\in \n$, 
then by taking limit over the ultrafilter $\U$ on $\n$ and applying Theorem \ref{thm-contraction-free}, 
Theorem \ref{thm-free-Haar-circular} and (1) of Proposition \ref{prop-ell-norm} we have 
\begin{align*}
\begin{split}
\norm{\sum^n_{i,j=1} \phi_{ij} \otimes ux_{ij}}_{L_2(\M;E)} & \leq \Phi_{unc}(E)
\norm{\sum^n_{i=1} \phi_{ij} \otimes \lambda(g_{ij}) \otimes ux_{ij}}_{L_2(\M \overline{\otimes} \N ;E)}\\
& \leq 2d_f(E)\norm{\sup_i \abs{\phi_i}}_2\Phi_{unc}(E)
\norm{\sum^n_{i=1} \lambda(g_{ij}) \otimes ux_{ij}}_{L_2(\N ;E)}\\
& \leq 2C_Gd_f(E)\norm{\sup_i \abs{\phi_i}}_2\Phi_{unc}(E)
\norm{\sum^n_{i=1} \WT_{ij} \otimes ux_{ij}}_{L_2(N ;E)}\\
& = 2C_Gd_f(E)\norm{\sup_i \abs{\phi_i}}_2\Phi_{unc}(E) \ell_{\WT}(uv)\\
& \leq 2C_Gd_f(E)\norm{\sup_i \abs{\phi_i}}_2\Phi_{unc}(E) \ell_{\WT}(u)\norm{v}\\
& = 2C_Gd_f(E)\norm{\sup_i \abs{\phi_i}}_2\Phi_{unc}(E) \ell_{\WT}(u) \norm{(x_{ij})}_{S^n_2 \otimes_{\min} OH_k},
\end{split}
\end{align*}
where $v : S^n_2 \rightarrow OH_k, \,\, e_{ij} \mapsto x_{ij}$.

Thus, by Proposition \ref{prop-replacing-2summing} we have
$$\pi^o_2(u) \leq 2C_Gd_f(E)^2\norm{\sup_i \abs{\phi_i}}_2\Phi_{unc}(E) \ell_{\WT}(u)$$ and consequently 
$${\WT C}^{oh}_2(E) \leq 2C_Gd_f(E)^2\norm{\sup_i \abs{\phi_i}}_2\Phi_{unc}(E).$$

\end{proof}

Finally we prove our main theorem. We state the theorem again with further comment on the constant.

\begin{thm}\label{thm-UncImpliesOH}
Let $\M$ and $\Phi=(\phi_i)_{i\geq 1}$ be the same as in Proposition \ref{prop-UncImpliesW-cotype2}.
Then, $E$ is $\Phi$-unconditional if and only if $E$ is completely isomorphic to an operator Hilbert space.
Moreover, we have $$\Phi_{unc}(E) \leq d_{cb}(E,OH(I)) 
\leq C_1 C_2^2 (1+\log C_2)$$ 
for some universal constant $C_1 > 0$ and $C_2 = d_f(E)^4\Phi_{unc}(E)\norm{\sup_i \abs{\phi_i}}_2$.
\end{thm}
\begin{proof}
Let $F$ be any finite dimensional subspace of $E$. Then $F$ is clearly  $\Phi$-unconditional with 
$\Phi_{unc}(F) \leq \Phi_{unc}(E)$. Thus, by Proposition \ref{prop-OH-TypeCotypeCharacterization} 
and \ref{prop-type-duality} we have
\begin{align*}
\begin{split}
d_{cb}(F,OH_{\text{dim}F}) & \leq d_f(F)\WT T^{oh}_2(F)\WT C^{oh}_2(F)\\
& \leq d_f(E)^2\WT K(F) \WT C^{oh}_2(F^*)\WT C^{oh}_2(F).
\end{split}
\end{align*}

$\WT C^{oh}_2(F)$ is estimated by Proposition \ref{prop-UncImpliesW-cotype2}, 
and a similar estimate for $\WT C^{oh}_2(F^*)$ can be done as follows. 

For fixed $n \in \n$, $(y_i)^n_{i=1} \subseteq F^*$ and unitaries $(u_i)^n_{i=1} \subseteq M_m$ 
we have by Proposition \ref{prop-duality} and the orthonormality of $\Phi$ that
\begin{align*}
\begin{split}
\lefteqn{\norm{\sum^n_{i=1}\phi_i \otimes y_i}_{L_2(\M;F^*)} \leq d_f(F)\norm{\sum^n_{i=1}\phi_i \otimes y_i}_{L_2(\M;E)^*}}\\
& = d_f(E)\sup\Big\{ \frac{\abs{\left\langle \sum^n_{i=1}\phi_i \otimes y_i, \sum^k_{j=1}\phi_j \otimes x_j \right\rangle}}
{\norm{\sum^k_{j=1}\phi_j \otimes x_j}_{L_2(\M;F)}} \Big\}\\
& \leq d_f(E)\Phi_{unc}(E)\sup\Big\{ \frac{\abs{\left\langle \sum^n_{i=1}\phi_i \otimes u_i \otimes y_i, 
\sum^k_{j=1}\phi_j \otimes u_j \otimes x_j \right\rangle}}
{\norm{\sum^k_{j=1}\phi_j\otimes u_j \otimes x_j}_{L_2(\M \overline{\otimes}M_m;F)}}\Big\}\\
& \leq d_f(E)\Phi_{unc}(E)\norm{\sum^n_{i=1}\phi_i \otimes u_i \otimes y_i}_{L_2(\M \overline{\otimes}M_m;F)^*}\\
& \leq d_f(E)\Phi_{unc}(E)\norm{\sum^n_{i=1}\phi_i \otimes u_i \otimes y_i}_{L_2(\M \overline{\otimes}M_m;F^*)}.
\end{split}
\end{align*}
By the same argument as in the proof of Proposition \ref{prop-UncImpliesW-cotype2} we get
$$\WT C^{oh}_2(F^*) \leq 2C_G d_f(E)^3\Phi_{unc}(E)\norm{\sup_i \abs{\phi_i}}_2.$$
Using this estimate and Proposition 12.4. of \cite{TJ89} we get 
\begin{align*}
\begin{split}
d_{cb}(F,OH_{\text{dim}F}) & \leq C \WT K(F) \leq C K(S_2(F)) \leq C_K C \log (1 + d(S_2(F),OH))\\
& = C_K C \log (1 + d_{cb}(F,OH_{\text{dim}F})),
\end{split}
\end{align*}
where $C = 4C^2_G d_f(E)^8\Phi_{unc}(E)^2\norm{\sup_i \abs{\phi_i}}^2_2$ and $C_K$ is the universal constant in Proposition 12.4. of \cite{TJ89}.
Thus, we have $d_{cb}(F,OH_{\text{dim}F}) \leq C_K C \log(2C_K C +1)$ for any finite dimensional subspace $F$ of $E$, and consequently
$$d_{cb}(E,OH(I)) \leq C_K C \log(2C_K C +1)$$ for some index set $I$.
\end{proof}

\section*{Acknowledgements}
The author would like to express his thanks to Prof. Marius Junge for his valuable and kind comments.

\bibliographystyle{amsplain}
\providecommand{\bysame}{\leavevmode\hbox
to3em{\hrulefill}\thinspace}

\end{document}